\def \D{{\cal D }}

\def \O{{\cal O }}

\def \R {\hat{R}}
\def \RR {\tilde{R}}
\def \PP {\tilde{\Phi}}
\def \Md {M^\dagger}

\def\Bbb#1{{#1\kern-.647em #1}}

\def \del{\partial}
\def \delb{\bar{\partial}}
\def \Db{\bar{D}}
\def \deltb{\bar{\delta}}

\def \etb {\bar{\eta}}
\def \zb {\bar{z}}
\def \xb {\bar{x}}
\def \xib {\bar{\xi}}

\def \be{\begin{equation}}
\def \eq{\end{equation}}
\def \bea{\begin{eqnarray}}
\def \eqa{\end{eqnarray}}

\def \ab {{\bar{a}}}
\def \bb {{\bar{b}}}
\def \cb {{\bar{c}}}
\def \db {{\bar{d}}}
\def \r {\rho}

\documentstyle[12pt]{article}

\begin{document}
\begin{titlepage}
\begin{center}
October 16, 95
        \hfill  LBL-\\
          \hfill    UCB-PTH-95/38 \\
\vskip .3in

{\large \bf Geometry of the Quantum Complex Projective Space $CP_q(N)$}
\footnote{This work was supported in part by the Director, Office of
Energy Research, Office of High Energy and Nuclear Physics, Division of
High Energy Physics of the U.S. Department of Energy under Contract
DE-AC03-76SF00098 and in part by the National Science Foundation under
grant PHY-90-21139.}
\vskip .5in
{Chong-Sun Chu, Pei-Ming Ho\footnote{email address:
cschu@physics.berkeley.edu,
pmho@physics.berkeley.edu}
and Bruno Zumino}
\vskip .2in
{\em
Department of Physics \\University of California \\ and\\
   Theoretical Physics Group\\
    Lawrence Berkeley Laboratory\\
      University of California\\
    Berkeley, CA 94720}
\end{center}

\vskip .2in

\begin{abstract}
The quantum deformation $CP_q(N)$ of complex projective space
is discussed. Many of the features present in the case of the quantum
sphere can be extended. The differential and integral calculus
is studied and $CP_q(N)$ appears as a quantum K\"{a}hler
manifold. The braiding of several copies of $CP_q(N)$ is
introduced and the anharmonic ratios of four collinear points
are shown to be invariant under quantum projective
transformations. They provide the building blocks of all
projective invariants. The Poisson limit is also described.
\end{abstract}
\end{titlepage}

%%%%%%%%%%%%%%%%%%%%%%%%%%%%%%%%%%%%%%%%%%%%%%%%%%%%%%%
\section{Introduction}
In a recent paper \cite{CHZ} the quantum sphere was described as
a complex quantum manifold. Then, in \cite{CHZ2}, the braiding
of several copies of the quantum sphere was introduced and quantum
anharmonic ratios (cross ratios) of four points on the sphere were
defined which are invariant
under the fractional transformation which describes the coaction
of the quantum group $SU_q(2)$ on the complex coordinates $z, \zb$
on the quantum sphere. In the present paper we extend the results
of \cite{CHZ} and \cite{CHZ2} to higher dimensions. In
Secs. 2 and 3 we define the quantum projective space $CP_q(N)$
in terms of both homogeneous and inhomogeneous complex coordinates and
we study the differential calculus on it. $CP_q(N)$ is shown to be the quantum
deformation of a K\"{a}hler manifold with the Fubini-Study metric.
In Sec. \ref{Poisson} we consider the Poisson limit.
In Sec. \ref{inte} we study the integration of functions on $CP_q(N)$ and give
explicit formulas for the integrals. Then, in Sec. \ref{BraidedCPN}
we introduce the braiding of several copies of $CP_q(N)$ and in
Sec. \ref{ProjGeom} we study the anharmonic ratio of four collinear
points in $CP_q(N)$. Just as in the classical case these anharmonic
ratios can be shown to be the building blocks of all invariants under
quantum projective transformations. For this reasons we have given
Sec. \ref{ProjGeom} the title ``Quantum Projective Geometry''.

All formulas and derivations of \cite{CHZ} and \cite{CHZ2} can be easily
modified, with a few changes of signs, to describe the quantum unit disk
and the coaction of quantum $SU_q(1,1)$ on it, as well as the corresponding
invariant anharmonic ratios. This provides a quantum deformation of
the Bolyai-Lobachevski\v{i} non-Euclidean geometry and of the differential
calculus on the Bolyai-Lobachevski\v{i} plane.
We shall not write here the modified equations appropriate for this case,
which can be guessed very easily, but
we would like to mention that the commutation relations between the
variables $z$ and $\zb$ for the unit disk are appropriate for a representation
of $z$ and $\zb$ as bounded operators in a Hilbert space.
This is to be contrasted with the case of the quantum sphere where $z$ and
$\zb$
must be unbounded operators. In a perfectly analogous way all formulas
and derivations of the present paper can be easily modified, with a few changes
of sign, to describe a quantum deformation of various higher dimensional
non-Euclidean geometries. Again we shall not do this explicitly here and leave
it as an exercise for the reader.

A different deformed algebra of functions on the Bolyai-Lobachevski\v{i} plane
has been considered in \cite{KL}.
The algebra of functions on complex projective space has been considered
by a number of authors, see for example \cite{arik}, \cite{sto} and \cite{mey}.
What we have shown here is that  a rich construction of differential geometry
and projective geometry can be carried out on this space.
It is not hard to extend most of the results of the present paper to the case
of
quantum  Grassmannian manifolds.

\section{$CP_q(N)$ as a Complex Manifold}

\subsection{Complex Quantum Space Covariant Under $S~U_q~(~N~+~1~)$}

For completeness, we list here the formulas we shall need to
construct the complex projective space.
Remember that the $SU_q(N+1)$ symmetry can be represented
\cite{CZ} on the
complex quantum space $C_q^{N+1}$ with
coordinates $x_i, \xb^i, i=0,1,...,N$ which satisfy the relations
\be \label{xx} x_i x_j =q^{-1} \RR^{kl}_{ij} x_k x_l , \eq
\be \label{xbx} \xb^i x_j=q (\RR^{-1})^{ik}_{jl} x_k \xb^l, \eq
\be \label{xbxb} \xb^i \xb^j=q^{-1} \RR^{ji}_{lk} \xb^k \xb^l. \eq
Here $q$ is a real number, $\RR^{kl}_{ij} $ is the $GL_q(N+1)$ $\R$-matrix
with indices running from $0$ to $N$,
and $\xb^i=x_i^*$ is the $*$-conjugate of $x_i$.
The Hermitian length
\be \label{length} L=x_i \xb^i \eq
is real and central.

Derivatives $D^i, \Db_i$ can be introduced (the usual symbols
$\del^a, \delb_b$ are reserved below for the derivatives on
$CP_q(N)$ ) which satisfy
\bea &D^i x_j =\delta^i_j +q \RR^{ik}_{jl} x_k D^l,
     &D^i \xb^j =q (\RR^{-1})^{ji}_{lk} \xb^k D^l,\\
     &\Db_i \xb^j =\delta^i_j +q^{-1} (\RR^{-1})^{lj}_{ki} \xb^k \Db_l,
     &\Db_i x_j =q^{-1} \PP^{lk}_{ji} x_k \Db_l
\eqa
and
\be D^i D^j=q^{-1} \RR^{ji}_{lk} D^k D^l, \eq
\be D^i \Db_j =q^{-1} \PP^{ki}_{lj} \Db_k D^l, \eq
\be \Db_i \Db_j =q^{-1} \RR^{kl}_{ij} \Db_k \Db_l. \eq
Here we have defined
\be \PP^{ij}_{kl}= \RR^{ji}_{lk} q^{2(i-l)} =\RR^{ji}_{lk} q^{2(j-k)} \eq
which satisfies
\be \PP^{ri}_{sj} (\RR^{-1})^{jk}_{il}=
    (\RR^{-1})^{ri}_{sj} \PP^{jk}_{il}= \delta^r_l \delta^k_s
\eq
and (sum over the index $k$)
\be \PP^{ik}_{jk}=\delta^i_j q^{2i+1}, \eq
\be \PP^{ki}_{kj}=\delta^i_j q^{2(N-i)+1}. \eq

There is a symmetry of this algebra:
\bea  &q\rightarrow q^{-1}, \label{s1} \\
      &x_i\rightarrow q^{-2i}\xb^i, \\
      &\xb^i\rightarrow x_i, \\
      &D^i\rightarrow q^{2i}\Db_i
\eqa
and
\be \Db_i\rightarrow D^i. \label{s4} \eq
Exchanging the barred and unbarred quantities in (\ref{s1}) - (\ref{s4})
we get another symmetry which is the inverse of this one.

Using the fact that $L$ commutes with $x_i, \xb^i$, a $*$-involution can be
defined for $D^i$
\be (D^i)^*=-q^{-2i'}L^n \Db_i L^{-n}, \eq
where
\be i'=N-i+1 \eq
for any real number $n$.
The $*$-involutions corresponding to different $n$'s are
related to one another by the symmetry of conjugation by $L$
\be a\rightarrow L^{m}aL^{-m}, \eq
where $a$ can be any function or derivative and
$m$ is the difference in the $n$'s.

The differentials $\xi_i=dx_i, \xib^i=(\xi_i)^*$ satisfy:
\be \label{xxi} x_i \xi_j =q \RR^{kl}_{ij} \xi_k x_l, \eq
\be \label{xbxi} \xb^i \xi_j =q (\RR^{-1})^{ik}_{jl} \xi_k \xb^l \eq
and
\be \xi_i \xi_j = -q \RR^{kl}_{ij} \xi_k \xi_l, \eq
\be \xib^i \xi_j =-q (\RR^{-1})^{ik}_{jl} \xi_k \xib^l. \eq

All the above relations are covariant under the transformation
\bea \label{x-transf}
      &x_i \rightarrow x_j T^j_i, &\xb^i \rightarrow (T^{-1})^i_j \xb^j,\\
      &D^i \rightarrow (T^{-1})^i_j D^j,
      &\Db_i \rightarrow  \Db_j q^{2i'} T^j_i q^{-2j'},\\
      &\xi_i \rightarrow \xi_j T^j_i, &\xib^i \rightarrow (T^{-1})^i_j \xib^j,
\eqa
where $T^i_j \in SU_q(N+1)$.

The exterior derivatives $\delta=\xi^i D_i, \deltb=\xib^i \Db_i$
on the holomorphic
and antiholomorphic functions
satisfy the undeformed Leibniz rule, $\delta^2=\deltb^2=0$ and
 $\deltb x_j =x_j \deltb$ etc.

\subsection{Projective Space $CP_q(N)$}

Define for $a=1,...,N$,
\footnote{The letters $a,b,c,e$ etc. run from 1 to $N$, while
$i, j, k, l$ run from 0 to $N$.}
\be z_a=x_0^{-1} x_a, \quad \zb^a=\xb^a (\xb^0)^{-1}. \eq
Since
\bea &x_0 x_a =q x_a x_0, & x_0 \xb^0 =\xb^0 x_0, \eqa
and \be x_0 \xb^a =q^{-1} \xb^a x_0, \eq
it follows from  (\ref{xx}) and (\ref{xbx}) that
\be
 \label{zz} z_a z_b =q^{-1}  \R^{ce}_{ab} z_c z_e, \eq
\be \label{zzb}
 \zb^a z_b=q^{-1} (\R^{-1})^{ac}_{be} z_c \zb^e -\lambda q^{-1} \delta^a_b.
\eq
where $\R^{ac}_{be}$ is the $GL_q(N)$ $\R$-matrix with indices running
from 1 to $N$ and $\lambda=q-1/q$.

Since
\bea &dz_a=x_0^{-1}(\xi_a-\xi_0 z_a), &d\zb^a=(\xib^a-\zb^a
\xib^0)(\xb^0)^{-1},
\eqa
and
\bea &x_0 \xi_0 =q^2 \xi_0 x_0, &x_0 \xib^0 =\xib^0 x_0, \eqa
it follows from (\ref{xxi}) and (\ref{xbxi}) that
\be z_a dz_b =q \R^{ce}_{ab} dz_c z_e \label{zdz}, \eq
\be \zb^a dz_b =q^{-1} (\R^{-1})^{ac}_{be} dz_c \zb^e, \label{zbdz} \eq
\be  dz_a dz_b = -q \R^{ce}_{ab} dz_c dz_e \label{dzdz} \eq
and
\be d\zb^a dz_b = -q^{-1} (\R^{-1})^{ac}_{be} dz_c d\zb^e. \label{dzdzb} \eq

The derivatives $\del^a, \delb_a$ are defined by requiring
$\delta \equiv dz_a \del^a$
and $\deltb \equiv d\zb^a\delb_a$ to be exterior differentiations.
%They can again be deduced from $D^i$ on the
%$SU_q(N+1)$ covariant quantum space as follows:
%\bea \label{delta-cpn}
%\delta &=& \sum_{i=0}^N \xi_i D^i \\ \nn
%  &=& \sum_{a=1}^N q dz_a x_0 D^a + \xi_0 D^0 + \xi_0 \sum_{a=1}^N z_a D^a \\
%%\nn
%  &=& \sum_{a=1}^N q dz_a x_0 D^a + q^{-1} \lambda^{-1} \xi_0 x_0^{-1} (\mu-1)
%\eqa
%where
%\be \mu=1+ q \lambda x_i D^i \eq
%is the scaling operator which satisfy
%\bea
%&\mu x_i =q^2 x_i \mu, &\mu \xb^i =\xb^i \mu ,\\
%&\mu D^i =q^{-2} D^i \mu, &\mu \Db_i = \Db_i \mu
%\eqa
%and so commutes with $z_a, \zb^a$.
%The second term in  (\ref{delta-cpn} can therefore be dropped on
%restricting to functions over $CP_q(N)$ and so we identify
%\be \label{del} \del^a=q x_0 D^a. \eq
%Similarly, $\deltb=d \zb^a \delb_a$ and
%\be \label{delb} \delb_a= \xb^0 \Db_a .\eq
%
It follows from (\ref{zdz}) and (\ref{zbdz}) that
\be \label{delz} \del^a z_b = \delta^a_b +q \R^{ac}_{be} z_c \del^e, \eq
\be \del^a \zb^b = q^{-1} (\R^{-1})^{ba}_{ec} \zb^c \del^e, \eq
\be \delb_a z_b =q \Phi^{ec}_{ba} z_c \delb_e, \eq
\be \delb_a \zb^b = \delta^b_a +q^{-1} (\R^{-1})^{eb}_{ca} \zb^c \delb_e, \eq
\be \label{deldel} \del^b \del^a =q^{-1} \R^{ab}_{ce} \del^e \del^c \eq
and
\be \del^a \delb_b=q \Phi^{ca}_{eb} \delb_c \del^e, \eq
where the $\Phi$ matrix is defined by
\be \Phi^{ca}_{db} =\R^{ac}_{bd} q^{2(c-b)}=\R^{ac}_{bd} q^{2(d-a)}. \eq

Similarly as in the case of quantum spaces
the algebra of the differential calculus on $CP_{q}(N)$
has the symmetry:
\bea   &q\rightarrow q^{-1}, \\
       &z_a\rightarrow q^{-2a}\zb^a, \\
       &\zb^a\rightarrow z_a, \\
       &\del^a\rightarrow q^{2a}\delb_a,
\eqa
and \be \delb_a\rightarrow \del^a. \eq

Also the $*$-involutions
\bea &z_a^*=\zb^a, \\
     &dz_a^*=d\zb^a,
\eqa
and
\be \label{inv-del} \del^{a*}=-q^{2n-2a'}\rho^{n} \delb_a \rho^{-n}, \eq
where
\be a'=N-a+1, \eq
and \be \r=1+\sum_{a=1}^{N}z_a\zb^a, \eq
can be defined for any $n$.
Corresponding to different $n$'s they are related with
one another by the symmetry of conjugation by $\rho$
to some powers followed by a recaling by appropriate powers of $q$.

In particular, the choice $n=N+1$ gives the $*$-involution
which has the correct classical limit of Hermitian conjugation
with the standard measure $\r^{-(N+1)}$ of $CP(N)$.

The transformation (\ref{x-transf}) induces a transformation on $CP_q(N)$
\be \label{z-transf} z_a \rightarrow
(T^0_0 +z_b T^b_0)^{-1} (T^0_a +z_b T^b_a). \eq
One can then calculate how the differentials transform
\be \label{dz-transf} dz_a \rightarrow dz_b M^b_a,
\quad d \zb^a \rightarrow (\Md)^a_b d \zb^b \eq
where $M^b_a$  is a matrix of function in $z_a$ with coefficients in
$SU_q(N+1)$
and $(\Md)^a_b \equiv (M^b_a)^*$.
Since $\delta, \deltb$ are invariant, it follows the transformation
on the derivatives
\be \label{del-transf} \del^a \rightarrow (M^{-1})^a_b \del^b, \quad
    (\del^a)^* \rightarrow (\del^b)^* ((\Md)^{-1})^b_a
\eq
The covariance of the $CP_q(N)$ relations under the transformation
(\ref{z-transf}), (\ref{dz-transf}) and (\ref{del-transf}) follows
directly from the covariance in $C_q^{N+1}$.

\section{A Note on the Differential Calculus}
In \cite{CHZ}, we showed that there exists a one form representation of the
differential.
The construction there can be generalized.
Let $A$ be a $*-$ involutive algebra  with coordinates $z_i, \zb_i$
 and differentials $dz_i, d \zb_i$ such that $\zb_i= z_i^*, d\zb_i=(dz_i)^*$.
If there exists a real element $a \in A$ and real unequal nonvanishing
constants $r, s$ such that
\be a z_i=r z_i a , \quad a dz_i =s dz_i a, \quad \forall i, \label{a}\eq
then, as easily seen,
\be \lambda \delta f =[\eta, f]_{\pm}, \quad
    \eta =  \frac{\lambda}{1-s/r} \delta a a^{-1}, \label{da}
\eq
\be \lambda \deltb f =[\etb, f]_{\pm}\quad
    \etb =  \frac{\lambda}{1-r/s} \deltb a a^{-1}, \label{dba}
\eq
and
\be \lambda d f =[\Xi, f]_{\pm}, \quad \Xi=\eta +\etb, \label{xi-eta} \eq
where $\pm$ applies for odd/even forms $f$.
Notice that (\ref{da}) and (\ref{dba}), and therefore (\ref{a}),
imply that
\be r a \delta a = s \delta a a, \quad r \deltb a a = s a \deltb a. \eq

It can be proved that $\eta^*=-\etb$ and so  $\Xi^*=-\Xi$. It  holds that
$\eta^2=\etb^2=0$. However $\Xi^2=\eta \etb + \etb \eta=
\lambda \delta \etb = \lambda \deltb \eta$  will generally
be nonzero.
Note that
\be \lambda d \Xi = [\Xi, \Xi]_+ =2 \Xi^2. \eq
Define
\be K= \delta \etb= \deltb \eta \label{K-eta} \eq
then
\be K= \frac{1}{2} d \Xi. \label{K-Xi} \eq
It follows that $dK=0$ and $K^*=K$. Thus in the case $K \neq 0$, we will call
it a K\"{a}hler form and $K^n$
\footnote{$n$ = complex dimension of the algebra. We consider only deformations
such that the Poincar\'{e} series of the deformed algebra and its classical
counterpart match.}
will be non-zero and define a real volume element for an
integral (invariant integral if $K^n$ is invariant).
$K$ also has the very nice property of commuting with
everything
\be K z_a =z_a K, \quad K dz_a =dz_a K. \eq
In the case of $S^2_q$, $K$ is just the area element.

Such a one-form representation for the calculus exists on both
$C_q^{N+1}$ and $CP_q(N)$.
For $C_q^{N+1}$,
we have
\bea &L x_i =x_i L, &L \xi_i = q^2 \xi_i L \eqa
and so by taking $a=L$, we have
\bea&\eta_0= -q^{-1}\delta L L^{-1},  &\etb_0 = q \deltb L L^{-1}. \eqa
In this case, $K$ is not the K\"{a}hler form one usually assigns to
$C_q^{N+1}$.
Rather, it gives $C_q^{N+1}$ the geometry of $CP_q(N)$ written in
homogeneous coordinates.

Similar relations hold for $CP_q(N)$ in inhomogeneous coordinates.
It is
\bea &\rho z_a =q^{-2} z_a \rho, &\rho dz_a =dz_a \rho \label{ro-zdz} \eqa
and therefore
\be \eta= -q^{-1}\delta \rho \rho^{-1}, \quad \etb= q \deltb \rho \rho^{-1}.
\eq
One can then compute
\bea K&=&\deltb \eta \\
    \label{metric}  &=&dz_a g^{a \bb} d \zb^b,
\eqa
where the metric $g^{a \bb}$ is
\be g^{a \bb} =q^{-1} \rho^{-2} (\rho \delta_{a b} - q^2 \zb^a z_b)
\label{gab}\eq
with inverse $g_{\bb c}$
\be g_{\bb c} g^{c \ab} = g^{a \cb} g_{\cb b} = \delta_{a b} \eq
given by
\be g_{\bb c} =q \rho (\delta_{b c}+ \zb^b z_c). \eq
This metric is the quantum deformation of the standard Fubini-Study metric
for $CP(N)$.

Notice that under the transformation (\ref{z-transf})
\be \eta \rightarrow \eta +q f^{-1} \delta f, \quad f=T^0_0 +z_b T^b_0 \eq
and so $K$ is invariant.
{}From (\ref{dz-transf}) and (\ref{metric}), it follows that
\be \label{g-transf1}
       {g^{a \bb} \rightarrow (M^{-1})^a_c  g^{c \db}  ((\Md)^{-1})^\db_\bb },
\eq
\be \label{g-transf2}
       {g_{\bb a} \rightarrow  (\Md)^\bb_\db  {g_{\db c}}  M^c_a}.
\eq

One can show that the volume element $dv_x$ in $C_q^{N+1}$
\bea dv_x &\equiv&\Pi_{j=0}^N (\xib^j L^{-1/2}) \Pi_{i=0}^N (L^{-1/2} \xi_i) \\
     &=&\rho^{-(N+1)} d\zb^N \cdots d \zb^1 dz_1 \cdots dz_N \cdot
     \xib^0 (\xb^0)^{-1} (x_0)^{-1} \xi_0.
\eqa
Since $dv_x$ is invariant, one can prove that
\be
   dv_z \equiv \rho^{-(N+1)} d\zb^N \cdots d \zb^1 dz_1 \cdots dz_N \label{dv}
\eq
is invariant also and is in fact equal to $K^N$ (up to a numerical factor).
The factor $\rho^{-(N+1)}$ justifies the choice $n=N+1$ for the
involution (\ref{inv-del}).

Having a quantum K\"{a}hler metric one can define connections, curvature,
a Ricci tensor and a Hodge star operation.
We shall not do it here because there seems to be no unique way to define
these constructs.
Still, once certian choices are made,
the full differential geometry can be developed.
See \cite{Ho} for a discussion of the quantum Riemannian case.

\section{Poisson Structures on $CP(N)$}\label{Poisson}

The commutation relations in the previous sections give us,
in the limit $q\rightarrow 1$,
a Poisson structure on $CP(N)$.
As usual, the Poisson Brackets (P.B.s) are obtained as the limit
\bea &(f,g)=\lim_{h\rightarrow 0}\frac{fg \mp gf}{h},
     &q=e^h=1+h+[h^2].
\eqa
It is straightforward to find
\bea   &(z_a, z_b)=z_a z_b, \hspace{1cm} a<b, \\
       &(z_a, \zb^b)= \left\{\begin{array}{ll}
                        z_a \zb^b, & a\neq b \\
                        2(1+\sum_{c=1}^{a}z_c \zb^c), &a=b
                      \end{array},\right.\\
       &(z_a, dz_b)=\left\{\begin{array}{ll}
                              z_a dz_b+2z_b dz_a, & a<b \\
                              2z_a dz_a, & a=b \\
                              z_a dz_b, & a>b
                          \end{array},\right. \\
       &(\zb^a, dz_b)=\left\{\begin{array}{ll}
                              -\zb^a dz_b, & a\neq b \\
                              -2\sum_{c=1}^{a}\zb^c dz_c, & a=b
                            \end{array}\right.
\eqa
and those following from the $*$-involution, which satisfies
\be (f, g)^* = (g^*, f^*). \eq

The P.B. of two differential forms $f$ and $g$
of degrees $m$ and $n$ respectively satisfies
\be (f, g) = (-1)^{mn+1} (g, f). \eq
The exterior derivatives $\delta, \deltb, d$  act on the P.B.s distributively,
for example
\be d(f, g) = (df, g) \pm (f, dg), \eq
where the plus (minus) sign applies for even (odd) $f$.
Notice that we have extended the concept of Poisson Bracket to include
differential forms.

In the classical limit (\ref{da}), (\ref{dba}) and (\ref{xi-eta}) become
\bea
   2\delta f = (\eta, f), \label{eta-f} \\
   2\deltb f = (\etb, f), \label{etb-f}
\eqa
and
\be
   2df = (\Xi, f). \label{xi-f}
\eq
Equations (\ref{K-eta}), (\ref{K-Xi}) and (\ref{ro-zdz}) to (\ref{dv})
are still valid, with $q=1$, but now
\be
   \Xi^2 = 0.
\eq
The Fubini-Study K\"{a}hler form
\be
   K=dz_a g^{a\bb}d\zb^b
\eq
has vanishing Poisson bracket with all functions and forms
and, naturally, it is closed.
We find the validity of (\ref{eta-f}), (\ref{etb-f}) and (\ref{xi-f})
very remarkable,
since the one-forms $\eta$, $\etb$ and $\Xi$ do not have to be
adjoined to the space of one-forms but already belong there naturally.

\section{Integration}\label{inte}

We now turn to the discussion of integration on $CP_q(N)$. We shall use the
notation $<f(z,\zb)>$ for
the right-invariant integral of a function $f(z,\zb)$ over $CP_q(N).$
It is defined, up to a normalization factor, by requiring
\be \label{inv-int} <\O f(z, \zb)>=0 \eq
for any left-invariant vector field $\O$ of $SU_q(N+1)$.
In \cite{CHZ}, the integral was computed for $CP_q(1)=S_q^2$
by considering explicitly how
the vector field act on functions. We shall follow a different and simpler
approach here.
First we notice that the identification
\be x_i /L^{1/2}=T^N_i, \quad \xb^i /L^{1/2}=(T^{-1})^i_N,
\quad i={0,1,...,N} \eq
reproduces (\ref{xx})-(\ref{length}). Thus if we define
\be \label{int} <f(z,\zb)> \equiv <f(z, \zb) |_{z_a=(T^N_0)^{-1} T^N_a,
    \zb^a=(T^{-1})^a_N /(T^{-1})^0_N}>_{SU_q(N+1)}, \eq
where $< \cdot >_{SU_q(N+1)}$ is the Haar measure \cite{Wor} on $SU_q(N+1)$,
then it follows immediately that (\ref{inv-int}) is satisfied.
\footnote{A similar startegy of using the ``angular'' measure to
define an integration has been employed by
H. Steinacker \cite{Ste} in constructing integration over the Euclidean space.}
Next we claim that
\be
\label{int1} <(z_1)^{i_1} ( \zb^1)^{j_1} \cdots (z_N)^{i_N} (\zb^N)^{j_N}> =0
{\;  \rm unless \;}
    i_1=j_1,..., i_N=j_N.
\eq
This is  because the integral is invariant under the finite transformation
(\ref{z-transf}).
For the particular choice $T^i_j =\delta^i_j \alpha_i$, with
$|\alpha_i|=1, \Pi_{i=0}^N \alpha_i =1$, this gives
\be z_a \rightarrow (\alpha_a /\alpha_0) z_a \eq
and so (\ref{int1}) follows.

In \cite{Wor}, Woronowicz proved the following interesting property
for the Haar measure
\be \label{int2} <f(T) g(T)>_{SU_q(N+1)}= < g(T) f(DTD) >_{SU_q(N+1)} \eq
where
\be (DTD)^i_j =D^i_k T^k_m D^m_j \eq
and
\be D^i_j = q^{-N+2i} \delta^i_j \eq
is the $D$ matrix for $SU_q(N+1)$.
It follows from (\ref{int2}) that
\be <f(z, \zb) g(z, \zb)> = <g(z, \zb) f(\D z, \D^{-1} \zb)> \eq
where
\be \D^a_b =\delta^a_b q^{2a}, \quad a,b =1,2,...,N .\eq

Introducing
\bea &\rho_r= 1+\sum_{a=1}^r z_a \zb^a, \eqa
one finds from (\ref{zz}) and (\ref{zzb}) that
\be \label{rhoz}  \rho_r z_a = \left\{ \begin{array}{ll}
                          z_a \rho_r & r<a \\
                          q^{-2}  z_a \rho_r & r \geq a
                          \end{array}
                 \right. ,
\eq
\be \rho_r \rho_s =\rho_s \rho_r \eq
and
\be \zb^a z_a= q^{-2} \rho_a - \rho_{a-1} \quad \mbox{(no sum).} \label{roro}
\eq

Because of (\ref{int1}), it is sufficient to determine integrals of the form
\be
   <{\rho_1}^{-i_1} \cdots {\rho_N}^{-i_N} >. \label{ro-int}
\eq
The values of the integers $i_a$ for (\ref{ro-int}) to make sense
will be determined later.

Consider
\bea <\zb_a {\rho_1}^{-i_1} \cdots {\rho_N}^{-i_N}  z_a >&=&
        <{\rho_1}^{-i_1} \cdots {\rho_N}^{-i_N}  z_a (q^{-2a}
\zb^a)>\nonumber\\
     &=& q^{-2a} <{\rho_1}^{-i_1} \cdots {\rho_N}^{-i_N}  (\rho_a-\rho_{a-1})>
\eqa
Using (\ref{rhoz}),
\bea \mbox{L.S. }
&=&q^{2(i_a+\cdots + i_N)} <{\rho_1}^{-i_1} \cdots {\rho_N}^{-i_N} \zb^a z_a >
                \nonumber\\
                &=& q^{2I_a} <{\rho_1}^{-i_1} \cdots {\rho_N}^{-i_N} \zb^a z_a
>
\eqa
where we have denoted
\be I_a =i_a + \cdots + i_N. \eq
Using (\ref{roro})
we get the recursion formula
\bea
&<{\rho_1}^{-i_1} \cdots {\rho_{a-1}}^{-i_{a-1}+1} {\rho_a}^{-i_a}
\cdots {\rho_N}^{-i_N} >[I_a+a] \nonumber \\
&=<{\rho_1}^{-i_1} \cdots {\rho_{a-1}}^{-i_{a-1}} {\rho_a}^{-i_a+1}
  \cdots {\rho_N}^{-i_N} >[I_a+a-1],
\eqa
where
\be [x] =\frac{q^{2x}-1}{q^2-1}. \eq
It is obvious then that
\be
    <{\rho_1}^{-i_1} \cdots {\rho_a}^{-i_a}  >
    =  <{\rho_1}^{-i_1} \cdots {\rho_{a-1}}^{-i_{a-1}-i_a} >
                                           \frac{[a]}{[I_a+a]}.
\eq
By repeated use of the recursion formula,
$<{\rho_1}^{-i_1} \cdots {\rho_N}^{-i_N} >$ reduces finally to
$<{\rho_1}^{-i_1 -i_2 \cdots  -i_N}>$
and
\be <{\rho_1}^{-I_1}> =  \frac{1}{[I_1+1]} <1>. \eq

Therefore
\be
   <{\rho_1}^{-i_1} \cdots {\rho_N}^{-i_N} >=<1>\Pi_{a=1}^N
\frac{[a]}{[I_a+a]}.
\eq
For this to be positive definite, $i_a$ should be restricted such that
$I_a +a > 0$ for $a=1,\cdots,N$.

\section{Braided $CP_{q}(N)$} \label{BraidedCPN}

The braiding of the $C_q^{N+1}$ quantum planes induces
a braiding on the $CP_{q}(N)$'s.
Let the first copy of quantum plane be denoted by $x_i, \xb^i$
and the second by $x'_i, \xb'^i$.

A consistent and covariant choice of commutation relations
between them is
\bea &x_i x'_j = \tau\RR_{ij}^{kl} x'_k x_l, \label{braid-xx}\\
     &\xb^i x'_j = \nu(\RR^{-1})_{jl}^{ik} x'_k \xb^l \label{braid-xxb}
\eqa
and their $*$-involutions for arbitrary numbers $\tau, \nu$.
If we choose $\tau=\nu^{-1}$ then the Hermitian length $L$
will remain central, $Lf'=f'L$, for any function $f'$ of $x', \xb'$.
%For the same reason as we have stated in \cite{CHZ2},
However, $L'$ does not commute with $x, \xb$.

Assuming that the exterior derivatives of the two copies
satisfy the Leibniz rule
\bea \delta'f=\pm f\delta',
        & \deltb'f=\pm f\deltb', \\
     \delta f'=\pm f'\delta,
        & \deltb f'=\pm f'\deltb,
\eqa
where the plus (minus) signs apply for even (odd) $f$ and $f'$,
and
\bea \delta\delta'=-\delta'\delta, & \delta\deltb'=-\deltb'\delta, \\
     \deltb\delta'=-\delta'\deltb, & \deltb\deltb'=-\deltb'\deltb,
\eqa
we obtain the commutation relations between functions and forms
of different copies by letting $\delta,\deltb,\delta'$ and $\deltb'$
to act on  (\ref{braid-xx}) and (\ref{braid-xxb}).
As usual, the commutation relations between derivatives and functions
of different copies can also be derived from
the commutation relations between differential forms and functions
using the Leibniz rule of the exterior derivatives and
the identifications $\delta=dx_i D^i, \deltb=d\xb^i\Db_i$
for both copies.

{}From the above we derive the braiding relations of
two braided copies of $CP_{q}(N)$ in terms of the inhomogeneous
coordinates. They are independent of the particular choice of
$\tau$ and $\nu$. We have
\bea   &z_a z'_b = q\R_{ab}^{ce} (z'_c -q^{-1}\lambda z_c)z_e,
        \label{zz-braid}\\
       &\zb'^a z_b = q^{-1}(\R^{-1})_{be}^{ac} z_c \zb'^e
        - q^{-1}\lambda\delta^a_b \label{zzb-braid}
\eqa
and their $*$-involutions as well as the commutation relations
between functions and forms of different copies
following the assumption that their exterior derivatives anticommute.

\section{Quantum Projective Geometry }\label{ProjGeom}

We will show in the following that many concepts of projective geometry
have an analogue in the deformed case, in particular we shall study
the deformed anharmonic ratios
(cross ratios).

\subsection{ Collinearity
 Condition }

Classically the collinearity conditions for $m$ distinct points can be given
in terms of the inhomogeneous coordinates
$\{z^A_a| A=1,2,\cdots,m; a=1,2,\cdots,N\}$ as:
\be
   (z^A_a-z^B_a)(z^C_a-z^D_a)^{-1}=(z^A_b-z^B_b)(z^C_b-z^D_b)^{-1},
   \label{ccc}
\eq
where  $A\neq B, C\neq D=1,\cdots,m$ and $a,b=1,\cdots,N$.

In the deformed case, the coordinates $\{z^A_a\}$ of $m$ points must
be braided for the commutation relations to be covariant, namely,
\be
   z^A_a z^B_b=q\R_{ab}^{ce}(z^B_c - q^{-1}\lambda z^A_c)z^A_e,\quad A\leq B,
   \label{ab-braid}
\eq
as an extension of (\ref{zz-braid}).
Equation (\ref{zzb-braid}) can also be generalized in the same way,
but we shall not need it in this section.
This braiding has the interesting property that
the algebra of $CP_{q}(N)$ is {\em self-braided}, that is,
(\ref{ab-braid}) allows the choice $A=B$.
This property makes it possible to talk about the coincidence of points.
Actually, the whole differential calculus for braided $CP_q(N)$
described in Sect.\ref{BraidedCPN} has this property.

Another interesting fact about this braiding is that for a fixed index
$a$ the commutation relation is identical to that for braided $S_q^2$
\footnote{This formula differs from the corresponding one in \cite{CHZ2}
because in the present paper we have used different ordering conventions.}:
\be
   z^A_a z^B_a=q^2 z^B_a z^A_a-q\lambda z^A_a z^A_a, \quad A\leq B.
   \label{sq2zz}
\eq

Since there is no algebraic way to say that two "points" are distinct
in the deformed case,
the collinearity conditions should avoid using expressions like
$(z^A_a-z^B_a)^{-1}$,
which are ill defined.
Denoting
\be [AB]_a := z^A_a - z^B_a, \eq
the collinearity conditions in the deformed case can be formulated as:
\be
   [AB]_a [CD]_b = q^2 [CD]_a [AB]_b, \quad \forall a,b,
   \label{collinear}
\eq
and $A<B\leq C<D$.
By (\ref{ab-braid}) this equation is formally equivalent to
the quantum counterpart of (\ref{ccc}):
\be
   [AB]_a [CD]_a^{-1} = [AB]_b [CD]_b^{-1},
   \label{collinear2}
\eq
where the ordering of $A,B,C,D$ is arbitrary.
The advantage of this formulation is that (\ref{collinear}) is a
quadratic polynomial condition and polynomials are well defined in the
braided algebra.

Therefore the algebra $Q$ of functions of $m$ collinear points
is the quotient of the algebra $A$ of $m$ braided copies of $CP_q(N)$
over the ideal $I:=\{f\alpha g|\forall f,g\in A; \forall \alpha\in CC\}$,
generated by $\alpha$ which stands for the collinearity conditions
(\ref{collinear}), i.e.,
$\alpha\in CC:=\{[AB]_a [CD]_b - q^2 [CD]_a [AB]_b| A<B\leq C<D\}$.

Two requirements have to be checked for this definition $Q:=A/I$ to make sense.
The first one is that for any $f\in A$ and $\alpha\in CC$,
\be
   f\alpha = \sum_i \alpha_i f_i, \forall f\in A, \label{f-alpha}
\eq
for some $f_i \in A$ and $\alpha_i\in CC$.
This condition ensures that the ideal $I$ generated by
the collinearity conditions is not "larger" than what we want,
as compared with the classical case.

Note that not all the collinearity conditions are independent.
In fact, it is sufficient (for formal manipulations, at least)
to consider only $B=C=m-1, D=m$ in either
(\ref{collinear}) or (\ref{collinear2}).
That is, we need only two points to fix a line.

We now check that (\ref{f-alpha}) is satisfied.
Obviously we only have to consider the cases $f=z^E_c$,
for arbitrary $E$ and $c$.
Let $\alpha(AB)_{ab}:=[AB]_a [CD]_b - q^2 [CD]_a [AB]_b$,
for $C=m-1$ and $D=m$.
Using (\ref{ab-braid}) one finds, after considerable algebra,
for $B\leq A<C<D$,
\be
   z^B_a\alpha(AC)_{bc}=q^2\R_{ab}^{he}\R_{ec}^{fg}
   \alpha(AC)_{hf}z^B_g.
\eq
For $A\leq B\leq C<D$, one finds similarly
\be
   z^B_a\alpha(AC)_{bc}=q^2\R_{ab}^{he}\R_{ec}^{fg}(\alpha(AC)_{hf}z^B_g
   +q^{-1}\lambda\alpha(AB)_{hf}[AB]_g).
\eq
Hence (\ref{f-alpha}) is proven for $B\leq C$.
Using
\be
   [CD]_a\alpha(AC)_{bc}=
(\R^{-1})_{ab}^{he}(\R^{-1})_{ec}^{fg}\alpha(AC)_{hf}[CD]_g,
\eq
for $B=D$ and
\be
  [BD]_a\alpha(AC)_{bc}=
  q^{-2}(\R^{-1})_{ab}^{he}(\R^{-1})_{ec}^{fg}\alpha(AC)_{hf}[BD]_g,
\eq
for $B>D$, together with the above two equations we immediately see that
(\ref{f-alpha}) is satisfied for $f=z^B_a$ also for $B\geq D$.
Therefore the first requirement is satisfied.

The second requirement is the invariance of $I$ under
the fractional transformations (\ref{z-transf}).
While this can be directly checked for (\ref{collinear}),
it is equivalent but simpler to consider another expression of the
collinearity conditions:
\be
   [AB]_a^{-1}[AB]_b = [CD]_a^{-1}[CD]_b,
\eq
where the ordering of $A,B,C,D$ is arbitrary.
Again we only have to consider the independent cases: $B<A=C=m-1$, $D=m$.
The fractional transformation has
\be
   [AB]_a\rightarrow -U(B)^{-1}[AB]_b z^A_c M_a^{bc}V(A)^{-1},\label{AB-transf}
\eq
where $U(B)=T^0_0+z^B_eT^e_0$, $V(A)=T^0_0+q z^A_f T^f_0$ and
$M_a^{bc}=T^b_0 T^c_a - q T^b_a T^c_0$.
So
\bea
   & & [AB]_a^{-1}[AB]_b \rightarrow\nonumber \\
   & & V(A)([AB]_c z^A_h M_a^{ch})^{-1}([AB]_e z^A_f M_b^{ef})V(A)^{-1}
       \nonumber\\
   &=& V(A)([AB]_c [AC]_c^{-1} [AC]_c z^A_h M_a^{ch})^{-1}
      ([AB]_e [AC]_e^{-1} [AC]_e z^A_f M_b^{ef})V(A)^{-1} \nonumber\\
   &=& V(A)([AC]_c z^A_h M_a^{ch})^{-1}([AB]_g [AC]_g^{-1})^{-1}
     ([AB]_s [AC]_s^{-1})([AC]_e z^A_f M_b^{ef})V(A)^{-1} \nonumber\\
   &=& V(A)([AC]_c z^A_h M_a^{ch})^{-1}([AC]_e z^A_f M_b^{ef})V(A)^{-1},
\eqa
(where we used (\ref{collinear2}) for the second equality)
which equals the transformation of $[AC]_a^{-1}[AC]_b$.
This means that the relation
$[AB]_a^{-1}[AB]_b - [AC]_a^{-1}[AC]_b = 0$
is preserved by the transformation.

\subsection{ Anharmonic Ratios }

Classically the anharmonic ratio of four collinear points is an invariant
of the projective mappings, which are the linear transformations
of the homogeneous coordinates.
In the deformed case, the homogeneous coordinates are the coordinates $x_i$,
$\bar{x}_i$ of the $SU_q(N+1)$-covariant quantum space,
and the linear transformations are the $GL_q(N+1)$ transformations
\footnote{The commutation relations (\ref{xx}) are also covariant under
$GL_q(N+1)$ transformations.}
(\ref{x-transf}), which induce the fractional transformations
(\ref{z-transf}) on the coordinates $z_i$,$\zb_i$ of the projective space
$CP_q(N)$.

We define the anharmonic ratio of $CP_q(N)$ for four collinear points
$\{z^A_a| A=1,2,3,4\}$ to be
\be
   [A1]_a [A4]_a^{-1} [B4]_a [B1]_a^{-1}, \label{cross-ratio}
\eq
where $A,B = 2,3$.
We wish to show that it is invariant.
Using (\ref{AB-transf}) and denoting $\tau(A):=[1A]_a [14]_a^{-1}$
which is independent of the index $a$ according to the collinearity condition,
we get
\be
   [AB]_a\rightarrow U(B)^{-1}(\tau(A)-\tau(B))P_a(A)V(A)^{-1},
\eq
where $P_a(A):=[14]_b z^A_c M_a^{bc}$.
Then the anharmonic ratio (\ref{cross-ratio}) transforms as
\bea
   [A1]_a [A4]_a^{-1} [B4]_a [B1]_a^{-1} &\rightarrow&
   U(1)^{-1}\tau(A)(1-\tau(A))^{-1}(1-\tau(B))\tau(B)^{-1}U(1)\nonumber\\
   &=& \tau(A)(1-\tau(A))^{-1}(1-\tau(B))\tau(B)^{-1},\nonumber\\
   &=& [A1]_a [A4]_a^{-1} [B4]_a [B1]_a^{-1},
\eqa
where we have used $z^1_a\tau(A) = \tau(A)z^1_a$ for any $A$,
which is true because we can represent $\tau(A)$ as $[1A]_a [14]_a^{-1}$
with the same index $a$ and then use $z^1_a [AB]_a = q^2 [AB]_a z^1_a$.

Because of the nice property (\ref{sq2zz}), we can use the results
about the anharmonic ratios of $S_q^2$ ( which is a special case of $CP_q(N)$
with $N=1$ but no collinear condition is needed there) in \cite{CHZ}.
Note that all the invariants as functions of $z^A_a$ for a fixed $a$
in $CP_q(N)$ are also invariants as functions of $z^A=z^A_a$ in $S_q^2$.
The reason is the following.
Consider the matrix $T^a_b$ defined by
\bea
   T^0_0=\alpha, & T^0_a=\beta, \\
   T^a_0=\gamma, & T^a_a=\delta,
\eqa
where $\alpha,\beta,\gamma,\delta$ are components of an $SU_q(2)$ matrix,
and $T^b_b=1$ for all $b\neq a$, both $>0$,
with all other components vanishing.
It is a $GL_q(N+1)$ matrix, but the transformation (\ref{z-transf})
of $z^A_a$ by this matrix is the fractional transformation on $S_q^2$
with coordinate $z^A=z^A_a$.

Therefore all the anharmonic ratios of $CP_q(N)$ must have a corresponding
anharmonic ratio of $S_q^2$.
On the other hand, since all the anharmonic ratios of $S_q^2$
are functions of only one of them \cite{CHZ2},
all them have a corresponding invariant of $CP_q(N)$,
which are functions of (\ref{cross-ratio}) and may also
be called anharmonic ratios.
Hence we have established a one to one correspondence between
the anharmonic ratios
of $S_q^2$ and $CP_q(N)$,
and as a consequence the fact that all the anharmonic ratios of $CP_q(N)$
are functions of only one of them.

The anharmonic ratios are important because they are the building blocks
of invariants in projective geometry.
For example, in the $n$-dimensional classical case for given $2(n+1)$ points
with homogeneous coordinates $\{x^A_i\}$, inhomogeneous coordinates $\{z^A_a\}$
where $A=1,\cdots,n$, $i=0,1,\cdots,n$, and $a=1,\cdots,n$,
we can construct an invariant
\be
   I:=\frac{(1,2,\cdots,n,n+1)(n+2,n+3,\cdots,2(n+1))}
           {(1,2,\cdots,n,n+2)(n+1,n+3,\cdots,2(n+1))}
   \label{I}
\eq
where $(A_0,\cdots,A_n)$ is the determinant of the matrix
\be
   \left(\begin{array}{ccc}
             1         & \cdots & 1         \\
             z^{A_0}_1 & \cdots & z^{A_n}_1 \\
             \vdots    & \ddots & \vdots    \\
             z^{A_0}_n & \cdots & z^{A_n}_n
         \end{array}\right).
\eq
It is invariant because $(A_0,\cdots,A_n)$ equals the determinant of the matrix
$M^i_j=x^{A_i}_j$, $i,j=0,\cdots,n$, divided by the factor
$x^{A_0}_0\cdots x^{A_n}_0$, which cancels between the numerator and
denominator
of $I$.
It can be shown that this invariant $I$ is in fact the anharmonic ratio
of four points $z,z',z^{n+1},z^{n+2}$, where $z$ ($z'$) is the intersection
of the line fixed by $z^{n+1},z^{n+2}$ with the $(n-1)$-dimensional subspace
fixed by $z^1,\cdots,z^n$ ($z^{n+3},\cdots,z^{2(n+1)}$).

It is remarkable that all this can also be done in the quantum case.
We can construct an invariant $I_q$ using the quantum determinant
and we can also formulate the condition for $(n+1)$ points to share an
$(n-1)$-dimensional subspace.
Furthermore, we know how to describe the intersection between subspaces
of arbitrary dimension spanned by given points.
It can be shown that the invariant $I_q$ is indeed an anharmonic ratio
in the same way as the classical case.

%%%%%%%%%%%%%%%%%%%%%%%%%%%%%%%%%%%%%%%%%%%%%%%%%%%%%%%%%%%%%%%%%%%%%%%

%%%%%%%%%%%%%%%%%%%%%%%%%%%%%%%%%%%%%%%%%%%%%%%%%%%%%%%%%%%%%%%%%%%%%%%

\end{document}